\input amstex
\documentstyle{amsppt}
\input bull-ppt
\keyedby{bull257/amh}

\topmatter
\cvol{26}
\cvolyear{1992}
\cmonth{Jan}
\cyear{1992}
\cvolno{1}
\cpgs{131-136}
\title Semistability of amalgamated products,\\
HNN-extensions, and all one-relator groups \endtitle
\author Michael L. Mihalik and Steven T. Tschantz\endauthor
\shortauthor{M. L. Mihalik and S. T. Tschantz}
\shorttitle{Amalgamated Products, HNN-Extensions, 
ONe-Relator Groups}
\address Department of Mathematics, Vanderbilt 
University, Nashville, Tennessee
37235\endaddress
\date August 29, 1990 and, in revised form, June 10, 
1991\enddate
\endtopmatter

\document
\heading 1. Introduction\endheading
\par Semistability at infinity is a geometric property 
used in the study of
ends of finitely presented groups. If a finitely 
presented group $G$ is
semistable at infinity, then sophisticated invariants for 
$G$, such as the
fundamental group at an end of $G$, can be defined (see 
\cite{10}). It is
unknown whether or not all finitely presented groups are 
semistable at
infinity, although by \cite{16} it suffices to know 
whether all 1-ended
finitely presented groups are semistable at infinity. 
There are a number of
results showing many 1-ended groups have this property, 
e.g., if $G$ is
finitely presented and contains a finitely generated, 
infinite, normal subgroup
of infinite index, then $G$ is semistable at infinity 
(see \cite{12} and for
other such results \cite{13--15}).
\par Semistability at infinity is of interest in the 
study of cohomology of
groups; if a finitely presented group $G$ is semistable 
at infinity, then
$H^2(G;\Bbb Z G)$ is free abelian (see \cite{6, 7}). This 
is conjectured to be
true for all finitely presented groups, but at present it 
is not even known for
2-dimensional duality groups (where one is discussing the 
dualizing module, see
\cite2).
\par For negatively curved groups (i.e., hyperbolic 
groups in the sense of
Gromov, see \cite8), semistability at infinity has 
additional interesting
consequences. If a negatively curved group $G$ is given 
the word metric with
respect to some finite generating set, then there is a 
compactification 
$\overline{G}$ of $G$ where a point of $\partial 
G=\overline{G}-G$ is a certain
equivalence class of proper sequences of points in $G$. 
The boundary of $G$ is
a compact, metrizable, finite-dimensional space, which 
determines the
cohomology of $G$. Bestvina and Mess have shown that if 
$G$ is a negatively
curved group, then for every ring $R$, there is an 
isomorphism of $RG$-modules
$H^i(G;RG)\cong\check H^{i-1}(\partial G; R)$ (\v Cech 
reduced). Geoghegan has
observed that results in \cite1 imply that a negatively 
curved group $G$ is
semistable at infinity iff $\partial G$ has the shape of 
a locally connected
continuum (see \cite6). Furthermore, in \cite1, ideas 
closely related to
semistability at infinity are
 used to analyze closed irreducible 3-manifolds with
negatively curved fundamental group.
\par A continuous map is {\it proper\/} if inverse images 
of compact sets are
compact. Proper rays $r$, $s\: [0,\infty)\rightarrow K$ 
in a locally finite
CW-complex $K$ are said to {\it converge to the same 
end\/} of $K$ if for every
compact $C\subseteq K$ there exists an $N$ such that 
$r([N,\infty))$ and
$s([N,\infty))$ are contained in the same path component 
of $K-C$. A locally
finite CW-complex $K$ is {\it semistable at infinity\/} 
if any two proper rays,
which converge to the same end of $K$, are properly 
homotopic. If $G$ is a
finitely presented group, then $G$ is {\it semistable at 
infinity\/} if for
some (equivalently any) finite CW-complex $X$ with 
$\pi_1(X)=G$, the universal
cover $\widetilde{X}$ of $X$ is semistable at infinity.
\par If $H$ is a subgroup of two groups $A$ and $B$, the 
amalgamated product
$A*_HB$ is the quotient of the free product of $A$ and 
$B$ where the copies of
$H$ in $A$ and $B$ are identified. If $H$ and $H'$ are 
isomorphic subgroups of
$A$, the HNN-extension $A*_H$ (where $H'$ is taken as 
given) is the quotient of
$A*\langle t\rangle $ where $H$ is identified with 
$tH't^{-1}$ (see \cite{11}).
In \cite{21}, Stallings proves a decomposition theorem 
for finitely generated
groups having more than one end in terms of amalgamated 
products or
HNN-extensions over finite subgroups. In \cite4, Dunwoody 
shows that for
finitely presented groups, the process of recursively 
applying this
decomposition theorem to the factor groups eventually 
terminates in 0-ended
(i.e., finite) and 1-ended factor groups. Our main result 
is the following:
\proclaim{Theorem 1} If $G=A*_HB$ is an amalgamated 
product where $A$ and $B$
are finitely presented and semistable at infinity, and 
$H$ is finitely
generated, then $G$ is semistable at infinity. If 
$G=A*_H$ is an HNN-extension
where $A$ is finitely presented and semistable at 
infinity, and $H$ is finitely 
generated, then $G$ is semistable at infinity.
\endproclaim
\par If $G$ is the fundamental group of a graph of groups 
(see \cite{20}), then
$G$ can be expressed as some combination of amalgamated 
products and
HNN-extensions of the vertex groups over the edge groups. 
Hence, if $G$ is the
fundamental group of a finite graph of groups in which 
each vertex group is
finitely presented and semistable at infinity and each 
edge group is finitely
generated, then $G$ is semistable at infinity. However, 
it is possible that a
group $G$ can be expressed as a combination of 
amalgamated products and
HNN-extensions of finitely presented groups over finitely 
generated (but not
finite) groups without $G$ being the fundamental group of 
a graph of groups
with these vertex and edge groups, hence the above 
theorem applies to a  larger
class of group decomposition. Although the question of 
semistability at
infinity for all finitely presented groups reduces to the 
same question for
1-ended groups, it is possible to obtain a 1-ended group 
$G=A*_HB$ where $A$,
$B$, and $H$ are infinite-ended (and similarly for 
HNN-extensions), and in fact
this is the essential difficulty in the proof of our main 
theorem.
\par As a corollary to the proof of Theorem 1, the same 
methods apply (with
homotopy replaced by homology in the sense of \cite7) to 
give a cohomology
version of this result.
\proclaim{Corollary 2} If $G=A*_HB$ is an amalgamated 
product where $A$ and $B$
are finitely presented, $H^2(A,\Bbb Z A)$ and $H^2(B;
\Bbb ZB)$ are free abelian, and
$H$ is finitely generated, then $H^2(G;\Bbb Z
 G)$ is free abelian. If $G=A*_H$ is
an HNN-extension where $A$ is finitely presented, 
$H^2(A;\Bbb Z A)$ is free
abelian, and $H$ is finitely generated, then $H^2(G;\Bbb 
Z G)$ is free abelian.
\endproclaim
\par As an application of our main result, we get the 
following general
theorem:
\proclaim{Theorem 3} All finitely generated one-relator 
groups are semistable
at infinity.
\endproclaim
\par Finally, as a corollary (using \cite 7 as before), 
we get a purely
cohomological result.

\proclaim{Corollary 4} If $G$ is a finitely generated 
one-relator group, then
$H^2(G;\Bbb Z G)$ is free abelian.
\endproclaim
\heading 2. Outline of proofs\endheading
\par We describe the proof of our main theorem in the 
amalgamated product case.
Take a presentation $P$ for $G=A\,*_HB$ by combining 
presentations for $A$ and
$B$, each containing generators for $H$. If $Z$ is the 
standard 2-complex
obtained from $P$, then $Z=X\cup Y$ where $X$ and $Y$ are 
subcomplexes of $Z$
with $\pi_1(X)=A$ and $\pi_1(Y)=B$, and $X\cap Y$ is a 
wedge of circles
representing generators for $H$ in both $\pi_1(X)$ and 
$\pi_1(Y)$. The
universal cover $\widetilde{Z}$ of $Z$ is a union of 
copies of $\widetilde{X}$
and $\widetilde{Y}$ attached along copies of the Cayley 
graph $\Gamma$ of $H$.
The group $G$ acts on the left of $\widetilde{Z}$, 
permuting copies of
$\widetilde{X}$, $\widetilde{Y}$, and $\Gamma$.
\par To prove the main theorem, we show that any two 
proper edge paths $r$ and
$s$ in $\widetilde{Z}$, converging to the same end of 
$\widetilde{Z}$, are
properly homotopic. The normal form structure of $A*_HB$ 
provides the geometric
structure to show that $r$ and $s$ are properly homotopic 
in case
$\operatorname{im}(r)\cup \operatorname{im}(s)$ 
intersects no copy of $\Gamma$
in an infinite set of vertices.
\par If $r$ or $s$ meets some copy of $\Gamma$, say 
$\Gamma_0$, in infinitely
many vertices, then (by replacing each ray with a 
properly homotopic ray passing
through these points) we may as well assume 
$V=\operatorname{im}(r)\cap
\operatorname{im}(s)\cap \Gamma_0$ contains infinitely 
many vertices. Let $q$
be a proper edge path in $\Gamma_0$ passing through 
infinitely many vertices in
$V$. Then $q$ and $r$ (and $s)$ converge to the same end 
of $\widetilde{Z}$,
and it suffices to show that $q$ and $r$ are properly 
homotopic (since 
then $q$ and $s$ are similarly properly homotopic, and 
thus $r$ and $s$ are
properly homotopic). Thus we are reduced to the case 
where one of our rays is
contained in a copy $\Gamma_0$ of $\Gamma$.

\par The main ideas in this, the main case in our work, 
are as follows. We
split $\widetilde{Z}$ into two connected pieces 
$\widetilde{Z}^+$ and
$\widetilde{Z}^-$, which intersect along $\Gamma_0$ by 
taking $\widetilde{X}_0$
and $\widetilde{Y}_0$ to be the copies of $\widetilde{X}$ 
and $\widetilde{Y}$
containing $\Gamma_0$ and then defining $\widetilde{Z}^+$ 
to be the component
of $(\widetilde{Z}-\widetilde{Y}_0)\cup\Gamma_0$ 
containing $\Gamma_0$, and $
\widetilde{Z}^-$ to be the component of 
$(\widetilde{Z}-\widetilde{X}_0)\cup
\Gamma_0$ containing $\Gamma_0$. By extracting 
ideas from the proof of Dunwoody's accessibility theorem 
\cite4, we show that a
certain configuration of rays and ends cannot occur in 
$\widetilde{Z}^+$ or 
$\widetilde{Z}^-$. (This configuration is represented in 
Figure 1, 
where $C$ is a compact set in $\widetilde{Z}$; $u$, $v$, 
$u'_i$, and $v'_i$ are
proper rays in $\Gamma_0$, with the $u'_i$ and $v'_i$ in 
different components of
$\Gamma_0-C$ and diverging from $u$ and $v$ at 
progressively later points, and
where ovals represent distinct ends of either 
$\widetilde{Z}^+$ or $
\widetilde{Z}^-$.) Because this configuration cannot 
occur, we can construct

\fighere{9pc}\caption{{\smc Figure\/} 1}

\noindent proper homotopies between any proper ray in 
$\Gamma_0$ and any proper ray in $
\widetilde{Z}^+$ or $(\widetilde{Z}^-)$ that converge to 
the same end of
$\widetilde{Z}^+$ (respectively, $\widetilde{Z}^-)$. In 
essence, this says that
the ends of $\widetilde{Z}^+$ and $\widetilde{Z}^-$, 
determined by $\Gamma_0$,
are semistable at infinity. This fact provides the 
geometric structure needed
to construct a patchwork of proper homotopies in 
$\widetilde{Z}$, giving a
proper homotopy between $r$ and $q$ and, thus, between 
the given $r$ and $s$.
\par The proof that all one-relator groups are semistable 
at infinity is by an
induction argument patterned after the proof by Magnus of 
the Freiheitssatz
(see \cite{11}). The proof makes use of our main theorem, 
the following
structure theorem for one-relator groups, and a simple 
fact about semistability
at infinity for factor groups in certain amalgamated 
products.
\proclaim{Lemma 5} Given any finitely generated one 
relator group $G$, there
exists a finite sequence of finitely generated one 
relator groups $H_1$,
$H_2,\dotsc, H_n=G$ such that, for each $i<n$, either 
$H_{i+1}$ or $H_{i+1}*
\Bbb Z$ is an HNN-extension of $H_i$ over a finitely 
generated group, and such
that $H_1$ is either a free group or else is isomorphic 
to a free product of a
free group  and a finite cyclic group.
\endproclaim
\proclaim{Lemma 6} If $G$ is finitely presented and 
$G*\Bbb Z$ is semistable at
infinity, then $G$ is semistable at infinity.
\endproclaim


\Refs
\ref\no 1 \by M. Bestvina and G. Mess \paper The boundary 
of negatively
curved groups \jour preprint, 1990\endref
\ref\no 2 \by R. Bieri \book Homological dimension of 
discrete groups
\publ Queen Mary College Math. Notes, London \yr 
1976\endref
\ref\no 3 \by W. Dicks and M. J. Dunwoody \book Groups 
acting on graphs
\publ Cambridge Univ. Press, Cambridge and New York \yr 
1989\endref
\ref\no 4 \by M. J. Dunwoody \paper The accessibility of 
finitely
presented groups \jour Invent. Math. \vol 81 \yr 1985 
\pages 449--457
\endref
\ref\no 5 \by H. Freudenthal \paper \"Uber die Enden 
topologischer
R\"aume und Gruppen \jour Math. Z. \vol 33 {\rm(1931)} 
\pages 692--713
\endref
\ref\no 6 \by R. Geoghegan \paper The shape of a group 
\inbook Geometric
and Algebraic Topology (H. Torunczyk, ed.) \publ Banach 
Center Publ.,
vol. 18, PWN, Warsaw \yr 1986 \pages 271--280\endref
\ref\no 7 \by R. Geoghegan and M. Mihalik \paper Free 
abelian cohomology
of groups and ends of universal covers \jour J. Pure 
Appl. Algebra
\vol 36 \yr 1985 \pages 123--137\endref
\ref\no 8 \by M. Gromov \paper Hyperbolic Groups \inbook 
Essays in
Group Theory (S. M. Gersten, ed.) \publ Math. Sci. Res. 
Inst. Publ. series,
vol. 8, Springer-Verlag, Berlin and New York \pages 
75--263\endref
\ref\no 9 \by H. Hopf \paper Enden offener Ra\"ume und 
unendliche
diskontinuierliche Gruppen \jour Comment. Math. Helv \vol 
16
\yr 1943 \pages 81--100\endref
\ref\no 10 \by B. Jackson \paper End invariants of groups 
extensions
\jour Topology \vol 21 \yr 1982 \pages 71--81\endref
\ref\no 11 \by R. C. Lyndon and P. E. Schupp \paper 
Combinatorial
group theory \jour Springer-Verlag, Berlin, Heidelberg 
and New York,
1970\endref
\ref\no 12 \by M. Mihalik \paper Semistability at the end 
of a group
extension \jour Trans. Amer. Math. Soc. \vol 277 \yr 1983
\pages 307--321\endref
\ref\no 13 \bysame \paper Ends of groups with the integers
as quotients \jour J. Pure Appl. Algebra \vol 35 \yr 1985
\pages 305--320\endref
\ref\no 14 \bysame \paper Ends of double extension groups
\jour Topology \vol 25 \yr 1986 \pages 45--53\endref
\ref\no 15 \bysame \paper Semistability at $\infty$ of 
finitely
generated groups and solvable groups \jour Topology Appl.
\vol 24 \yr 1986 \pages 259--269\endref
\ref\no 16 \bysame \paper Semistability at 
$\infty,\infty$-ended
groups and group cohomology \jour Trans. Amer. Math. Soc.
\vol 303 \yr 1987 \pages 479--485\endref
\ref\no 17 \by M. Mihalik and S. Tschantz \book
Semistability of amalgamated products and HNN-extensions
\publ Mem. Amer. Math. Soc. (to appear)\endref
\ref\no 18 \bysame \paper All one relator groups are 
semistable
at infinity \jour preprint, 1991\endref
\ref\no 19 \by G. P. Scott and T. Wall \paper Topological
methods in group theory \inbook London Math. Soc. Lecture
Note Ser., vol. 36 \publ Cambridge Univ. Press \yr 1979
\pages 137--203\endref
\ref\no 20 \by J.-P. Serre \book Trees \publ 
Springer-Verlag,
Berlin and New York \yr 1980\endref
\ref\no 21 \by J. Stallings \book Group theory and three 
dimensional
manifolds \publ Yale Math. Mono., vol. 4, Yale Univ. Press,
New Haven, CT \yr 1972\endref

\endRefs
\enddocument